\documentclass[12pt]{amsart}
\usepackage[cp850]{inputenc}
\usepackage{graphics}

\newtheorem{theorem}{Theorem}[section]
\newtheorem{lemma}[theorem]{Lemma}

\newtheorem{corollary}[theorem]{Corollary}

\theoremstyle{definition}

\newtheorem{example}[theorem]{Example}

\theoremstyle{remark}
\newtheorem{remark}[theorem]{Remark}

\usepackage{amssymb,amsmath}

\newcommand{\Ric}{\mbox{$\mathrm{Ric}$}}
\newcommand{\m}{\mbox{$M$}}
\newcommand{\R}[1]{\mbox{${\mathbb R}^{#1}$}}
\newcommand{\h}[1]{\mbox{${\mathbb H}^{#1}$}}
\newcommand{\s}[1]{\mbox{${\mathbb S}^{#1}$}}
\newcommand{\g}[2]{\mbox{$\langle #1 ,#2 \rangle$}}
\newcommand{\fle}{\mbox{$\rightarrow$}}

\newcommand{\rf}[1]{\mbox{(\ref{#1})}}
\newcommand{\rl}[1]{{~\ref{#1}}}
\newcommand{\nablao}{\mbox{$\nabla^{\mathrm{o}}$}}
\newcommand{\nablabar}{\mbox{$\overline{\nabla}$}}
\newcommand{\xm}{\mbox{$\mathcal{X}(\m)$}}
\newcommand{\bin}[1]{\mbox{$\binom{n}{#1}$}}

\def\beq{\begin{equation}}
\def\eeq{\end{equation}}

\newcommand{\Lx}{\mbox{$L_kx=Ax+b$}}
\newcommand{\Lxb}{\mbox{$L_kx=Ax$}}

\def\Hn{\mbox{$\mathbb H^{n+1}$}}
\def\Sn{\mbox{$\mathbb S^{n+1}$}}

\newcommand{\mbar}{\mbox{$\mathbb{M}^{n+1}_c$}}
\newcommand{\x}{\mbox{$x:M^n\rightarrow\mbar\subset\mathbb{R}^{n+2}_q$}}

\begin{document}
\title[Hypersurfaces in space forms satisfying $L_kx=Ax+b$]
{Hypersurfaces in space forms satisfying the condition $L_kx=Ax+b$}

%    Information for first author
\author{Luis J. Al\'\i as}
\address{Departamento de Matem\'{a}ticas, Universidad de Murcia, E-30100 Espinardo, Murcia, Spain}
\email{ljalias@um.es}
\thanks{L.J. Al\'\i as was partially supported by MEC project MTM2007-64504, and Fundaci\'{o}n S\'{e}neca project
04540/GERM/06, Spain. This research is a result of the activity developed within the framework of the Programme in Support of Excellence Groups of the Regi\'{o}n de Murcia, Spain, by Fundaci\'{o}n S\'{e}neca, Regional Agency
for Science and Technology (Regional Plan for Science and Technology 2007-2010).}

%    Information for second author
\author{S. M. B. Kashani}
\address{Department of Mathematics, Faculty of Sciences, Tarbiat Modares University, P.O. Box 14115-175, Tehran, Iran}
\email{kashanism@yahoo.com,kashanim@modares.ac.ir}
\thanks{}

%    General info

%\date{}

%\dedicatory{}

%\keywords{}

\begin{abstract}
We study hypersurfaces either in the sphere \s{n+1} or in the hyperbolic space \h{n+1} whose position vector $x$
satisfies the condition \Lx, where $L_k$ is the linearized operator of the
$(k+1)$-th mean curvature of the hypersurface for a fixed $k=0,\ldots,n-1$,
$A\in\R{(n+2)\times (n+2)}$ is a constant matrix and $b\in\R{n+2}$ is a constant
vector. For every $k$, we prove that when $A$ is self-adjoint and $b=0$, the only hypersurfaces satisfying that
condition are hypersurfaces with zero $(k+1)$-th mean curvature and constant $k$-th mean curvature, and open pieces of
standard Riemannian products of the form
$\s{m}(\sqrt{1-r^2})\times\s{n-m}(r)\subset\s{n+1}$, with $0<r<1$, and
$\h{m}(-\sqrt{1+r^2})\times\s{n-m}(r)\subset\h{n+1}$, with $r>0$. If $H_k$ is constant, we also obtain a classification
result for the case where $b\neq 0$.
\end{abstract}

\maketitle

\section{Introduction}
\label{s1}
In \cite{AG} and inspired by Garay's extension of Takahashi theorem \cite{Ta,Ch1,Ch2} and its subsequent generalizations and
extensions \cite{CP,Ga2,DPV,HV2,AFL1,AFL2}, the first author jointly with G\"urbuz started the study of hypersurfaces in the Euclidean
space satisfying the general condition $L_k x=Ax+b$, where $A\in\R{(n+1)\times (n+1)}$ is a constant matrix and
$b\in\R{n+1}$ is a constant vector (we refer the reader to the Introduction of \cite{AG} for further details). In
particular, the following classification result was given in \cite[Theorem 1]{AG}.
\begin{theorem}
Let $x:M^n\rightarrow\R{n+1}$ be an orientable hypersurface immersed into the Euclidean space
and let $L_k$ be the linearized operator of the $(k+1)$-th mean curvature of $M$, for
some fixed $k=0,\ldots,n-1$. Then the immersion satisfies the condition $L_kx=Ax+b$ for
some constant matrix $A\in\R{(n+1)\times(n+1)}$ and some constant vector
$b\in\R{n+1}$ if and only if it is one of the following hypersurfaces in \R{n+1}:
\begin{enumerate}
\item a hypersurface with zero $(k+1)$-th mean curvature,
\item an open piece of a round hypersphere $\s{n}(r)$,
\item an open piece of a generalized  right spherical cylinder $\s{m}(r)\times\R{n-m}$, with $k+1\leq m \leq n-1$.
\end{enumerate}
\end{theorem}

In this paper, and as a natural continuation of the study started in \cite{AG}, we
consider the study of hypersurfaces $M^n$ immersed either into the sphere
$\s{n+1}\subset\mathbb{R}^{n+2}$ or into the hyperbolic space
$\h{n+1}\subset\mathbb{R}^{n+2}_1$ whose position vector $x$ satisfies the condition
\Lx. Here and for a fixed integer $k=0,\ldots,n-1$, $L_k$ stands for the linearized operator of the $(k+1)$-th mean
curvature of the hypersurface, denoted by $H_{k+1}$, $A\in\R{(n+2)\times (n+2)}$ is a constant
matrix and $b\in\R{n+2}$ is a constant vector. For the sake of simplifying the
notation and unifying the statements of our main results, let us denote by \mbar\
either the sphere $\Sn\subset\mathbb{R}^{n+2}$ if $c=1$, or the hyperbolic space
$\Hn\subset\mathbb{R}^{n+2}_1$ if $c=-1$. In this new situation, the codimension of
the manifold $M^n$ in the (pseudo)-Euclidean space $\mathbb{R}^{n+2}_q$ where it is lying is 2, which
increases the difficulty of the problem. In the case where $A$ is self-adjoint and
$b=0$ we are able to give the following classification result.
\begin{theorem}
\label{th1}
Let \x\ be an orientable hypersurface immersed either into the Euclidean sphere
$\Sn\subset\mathbb{R}^{n+2}$ (if $c=1$) or into the hyperbolic space
$\Hn\subset\mathbb{R}^{n+2}_1$ (if $c=-1$), and let $L_k$ be the linearized
operator of the $(k+1)$-th mean curvature of \m, for some fixed $k=0,\ldots,n-1$. Then the immersion satisfies the
condition \Lxb\ for some self-adjoint constant matrix $A\in\R{(n+2)\times(n+2)}$ if and only if it is one of the following
hypersurfaces:
\begin{enumerate}
\item a hypersurface having zero $(k+1)$-th mean curvature and constant $k$-th mean curvature;
\item an open piece of a standard Riemannian product $\s{m}(\sqrt{1-r^2})\times\s{n-m}(r)\subset\s{n+1}$, $0<r<1$,
if $c=1$;
\item an open piece of a standard Riemannian product $\h{m}(-\sqrt{1+r^2})\times\s{n-m}(r)\subset\h{n+1}$, $r>0$,
if $c=-1$.
\end{enumerate}
\end{theorem}
Let us recall that every compact hypersurface immersed into the hyperbolic space $\mathbb{H}^{n+1}$ has an elliptic point, that is, a point where all the
principal curvatures are positive (for a proof see, for instance, \cite[Lemma 8]{AKS}). The same happens for every
compact hypersurface immersed into an open hemisphere $\mathbb{S}^{n+1}_{+}$  (see, for instance, \cite[Section 3]{AAR} for a proof in the case $n=2$, although the
proof works also in the general $n$-dimensional case). In particular, this implies that there exists no compact hypersurface either in $\mathbb{H}^{n+1}$ or in $\mathbb{S}^{n+1}_{+}$ with vanishing $(k+1)$-th mean curvature, for
every $k=0,\ldots,n-1$. Since the standard Riemannian products $\s{m}(\sqrt{1-r^2})\times\s{n-m}(r)\subset\s{n+1}$ are not contained in an open hemisphere, then we have the following non-existence result as a consequence of our
Theorem\rl{th1}.
\begin{corollary}
There exists no compact orientable hypersurface either in $\mathbb{H}^{n+1}$ or in $\mathbb{S}^{n+1}_{+}$ satisfying the condition
\Lxb\ for some self-adjoint constant matrix $A\in\R{(n+2)\times(n+2)}$, where $L_k$
stands for any of the linearized operators of the higher order mean curvatures.
\end{corollary}

When $k=1$ the operator $L_1$ is the
operator $\Box$ introduced by Cheng and Yau in \cite{CY} for the study of
hypersurfaces with constant scalar curvature. In that case, since the scalar curvature of \m\ is given
by $n(n-1)(c+H_2)$ (see equation \rf{scalar}) we get the following consequence.
\begin{corollary}
Let \x\ be an orientable hypersurface immersed either into the Euclidean sphere
$\Sn\subset\mathbb{R}^{n+2}$ (if $c=1$) or into the hyperbolic space
$\Hn\subset\mathbb{R}^{n+2}_1$ (if $c=-1$), and let $\Box$ be the Cheng and Yau operator on \m. Then the immersion
satisfies the
condition $\Box x=Ax$ for some self-adjoint constant matrix $A\in\R{(n+2)\times(n+2)}$ if and only if it is one of the
following hypersurfaces:
\begin{enumerate}
\item a hypersurface having constant scalar curvature $n(n-1)c$ and constant mean curvature;
\item an open piece of a standard Riemannian product $\s{m}(\sqrt{1-r^2})\times\s{n-m}(r)\subset\s{n+1}$, $0<r<1$,
if $c=1$;
\item an open piece of a standard Riemannian product $\h{m}(-\sqrt{1+r^2})\times\s{n-m}(r)\subset\h{n+1}$, $r>0$,
if $c=-1$.
\end{enumerate}
\end{corollary}

In particular, when $n=2$, and taking into account that the only surfaces either in \s{3} or \h{3} having constant
mean curvature and constant Gaussian (or scalar) curvature equal to the Gaussian curvature of the ambient space are
the totally geodesic ones, we obtain the following result.
\begin{corollary}
\label{co1}
Let $x:M^2\rightarrow\mathbb{M}^{3}_c\subset\mathbb{R}^{4}_q$ be an orientable surface immersed either into the Euclidean
sphere $\s{3}\subset\mathbb{R}^{4}$ (if $c=1$) or into the hyperbolic space
$\h{3}\subset\mathbb{R}^{4}_1$ (if $c=-1$), and let
$L_1=\Box$ be the Cheng and Yau operator of $M$. Then the immersion satisfies the
condition $\Box x=Ax$ for some self-adjoint constant matrix $A\in\R{(4)\times(4)}$ if and only if it is one of the
following surfaces:
\begin{enumerate}
\item an open piece of either a totally geodesic round sphere $\s{2}\subset\s{3}$ or a standard
Riemannian product $\s{1}(\sqrt{1-r^2})\times\s{1}(r)\subset\s{3}$, $0<r<1$, if $c=1$;
\item an open piece of either a totally geodesic hyperbolic plane  $\h{2}\subset\h{3}$ or a standard
Riemannian product $\h{1}(-\sqrt{1+r^2})\times\s{1}(r)\subset\h{3}$, $r>0$, if $c=-1$.
\end{enumerate}
\end{corollary}

\begin{remark}
A different but related result to our Theorem\rl{th1} has been proved recently by Yang and Liu in \cite{YL}.
In fact, instead of assuming that $A$ is self-adjoint, they assume that $H_k$ is constant and reach the same
classification. Specifically, they use the method of moving frames to derive the basic equations for the
hypersurface and then, following the techniques introduced by Al\'\i as, Ferrßndez and Lucas in \cite{AFL2}
for the case $k=0$ and extended by Al\'\i as and G\"urb\"uz in \cite{AG} for general $k$, they prove that the
hypersurface must be one of the standard examples.
\end{remark}

On the other hand, in the case where $A$ is self-adjoint and $b\neq 0$ we are able to prove the following classification
result.
\begin{theorem}
\label{th2}
Let \x\ be an orientable hypersurface immersed either into the Euclidean sphere
$\Sn\subset\mathbb{R}^{n+2}$ (if $c=1$) or into the hyperbolic space
$\Hn\subset\mathbb{R}^{n+2}_1$ (if $c=-1$), and let $L_k$ be the linearized
operator of the $(k+1)$-th mean curvature of \m, for some fixed $k=0,\ldots,n-1$. Assume that $H_k$ is constant.
Then the immersion satisfies the
condition \Lx\ for some self-adjoint constant matrix $A\in\R{(n+2)\times(n+2)}$ and some non-zero constant vector
$b\in\R{n+2}$ if and only if:
\begin{itemize}
\item[(i)] $c=1$ and it is an open piece of a totally umbilical round sphere $\s{n}(r)\subset\s{n+1}$, $0<r<1$.
\item[(ii)] $c=-1$ and it is one of the following hypersurfaces in \h{n+1}:
\begin{enumerate}
\item an open piece of a totally umbilical hyperbolic space $\h{n}(-r)$, $r>1$,
\item an open piece of a totally umbilical round sphere $\s{n}(r)$, $r>0$,
\item an open piece of a totally umbilical Euclidean space $\R{n}$.
\end{enumerate}
\end{itemize}
\end{theorem}

\section{Preliminaries}
\label{s2}
Throughout this paper we will consider both the case of hypersurfaces immersed into the Euclidean sphere
\[
\s{n+1}=\{x=(x_0,\dots,x_{n+1})\in\mathbb{R}^{n+2} : \g{x}{x}=1 \},
\]
and the case of hypersurfaces immersed into the hyperbolic space \h{n+1}. In this last case, it will be appropriate to
use the Minkowski space model of hyperbolic space. Write $\mathbb{R}^{n+2}_1$ for $\mathbb{R}^{n+2}$, with coordinates
$(x_0,\ldots,x_{n+1})$, endowed with the Lorentzian metric
\[
\g{}{}=-dx_0^2+dx_1^2+\cdots+dx_{n+1}^2.
\]
Then
\[
\h{n+1}=\{x\in\mathbb{R}^{n+2}_1 : \g{x}{x}=-1, x_0>0 \}
\]
is a complete spacelike hypersurface in $\mathbb{R}^{n+2}_1$ with constant sectional curvature $-1$ which provides
the Minkowski space model for the hyperbolic space.

In order to simplify our notation, we will denote by \mbar\ either the sphere
$\Sn\subset\mathbb{R}^{n+2}$ if $c=1$, or the hyperbolic space
$\Hn\subset\mathbb{R}^{n+2}_1$ if $c=-1$.  We will also denote by $\g{}{}$, without
distinction, both the Euclidean metric on $\mathbb{R}^{n+2}$ and the Lorentzian
metric on $\mathbb{R}^{n+2}_1$, as well as the corresponding (Riemannian) metrics
induced on \mbar\ and on \m. Consider \x\ (with $q=0$ if $c=1$, and $q=1$ if $c=-1$)
a connected orientable hypersurface immersed into \mbar\ with Gauss map $N$.
Throughout this paper we will denote by $\nablao$, $\nablabar$ and $\nabla$ the
Levi-Civita connections on \R{n+2}, \mbar\ and \m, respectively. Then, the basic
Gauss and Weingarten formulae of the hypersurface are written as
\[
\nablao_XY=\nablabar_XY-c\g{X}{Y}x= \nabla_XY+\g{SX}{Y}N-c\g{X}{Y}x
\]
and
\[
SX=-\nablabar_XN=-\nablao_XN,
\]
for all tangent vector fields $X,Y\in\xm$, where $S:\xm\rightarrow\xm$ stands for the
shape operator (or Weingarten endomorphism) of \m\ with respect to the chosen
orientation $N$. As is well known, $S$ defines a self-adjoint linear operator on
each tangent plane $T_p\m$, and its eigenvalues $\kappa_1(p), \ldots, \kappa_n(p)$
are the principal curvatures of the hypersurface. Associated to the shape operator
there are $n$ algebraic invariants given by
\[
s_k(p)=\sigma_k(\kappa_1(p), \ldots, \kappa_n(p)), \quad 1\leq k\leq n,
\]
where $\sigma_k:\mathbb R^n\to\mathbb R$ is the elementary symmetric
function in $\mathbb R^n$ given by
\[
\sigma_k(x_1,\ldots, x_n)=\sum_{i_1<\cdots<i_k}x_{i_1}\ldots x_{i_k}.
\]
Observe that the characteristic polynomial of $S$ can be writen in terms
of the $s_k$'s as
\beq
\label{poly}
Q_S(t)=\det(tI-S)=\sum_{k=0}^n(-1)^ks_kt^{n-k},
\eeq
where $s_0=1$ by definition. The $k$-th mean curvature $H_k$ of the hypersurface is then defined by
\[
\bin{k}H_k=s_k, \quad 0\leq k\leq n.
\]

In particular, when $k=1$ $H_1=(1/n)\sum_{i=1}^n\kappa_i=(1/n)\mathrm{trace}(S)=H$
is nothing but the mean curvature of \m, which is the main extrinsic curvature of
the hypersurface. On the other hand, $H_2$ defines a geometric quantity which is
related to the (intrinsic) scalar curvature of \m. Indeed, it follows from the Gauss
equation of \m\ that its Ricci curvature is given by
\[
\Ric(X,Y)=(n-1)c\g{X}{Y}+nH\g{SX}{Y}-\g{SX}{SY}, \quad X,Y\in\xm,
\]
and then the scalar curvature of \m\ is
\begin{eqnarray}
\label{scalar}
\mathrm{tr}(\Ric) & = & n(n-1)c+n^2H^2-\mathrm{tr}(S^2) \\
\nonumber {} & = &
n(n-1)c+\left(\sum_{i=1}^n\kappa_i\right)^2-\sum_{i=1}^n\kappa_i^2=n(n-1)(c+H_2).
\end{eqnarray}
In general, when $k$ is odd the curvature $H_k$ is extrinsic (and its sign
depends on the chosen orientation), while when $k$ is even the curvature $H_k$ is
intrinsic and its value does not depend on the chosen orientation.

The classical Newton transformations $P_k:\xm\rightarrow\xm$ are defined
inductively from the shape operator $S$ by
\[
P_0=I \quad \mathrm{and} \quad P_k=s_kI-S\circ P_{k-1}=\bin{k}H_kI-S\circ P_{k-1},
\]
for every $k=1\ldots,n$, where $I$ denotes the identity in \xm.
Equivalently,
\[
P_k=\sum_{j=0}^k(-1)^{j}s_{k-j}S^{j}=\sum_{j=0}^k(-1)^{j}\bin{k-j}H_{k-j}S^{j}.
\]
Note that by the Cayley-Hamilton theorem, we have $P_n=0$ from \rf{poly}. Observe
also that when $k$ is even, the definition of $P_k$ does not depend on the chosen
orientation, but when $k$ is odd there is a change of sign in the definition of
$P_k$.

Let us recall that each $P_k(p)$ is also a self-adjoint linear operator on
each tangent plane $T_p\m$ which commutes with $S(p)$. Indeed, $S(p)$ and
$P_k(p)$ can be simultaneously diagonalized: if $\{ e_1, \ldots, e_n\}$
are the eigenvectors of $S(p)$ corresponding to the eigenvalues
$\kappa_1(p), \ldots, \kappa_n(p)$, respectively, then they are also the
eigenvectors of $P_k(p)$ with corresponding eigenvalues given by
\beq
\label{muik}
\mu_{i,k}(p)=\frac{\partial \sigma_{k+1}}{\partial x_i}
(\kappa_1(p), \ldots, \kappa_n(p))=\sum_{i_1<\cdots<i_k,i_j\neq i}\kappa_{i_1}(p)\cdots\kappa_{i_k}(p),
\eeq
for every $1\leq i\leq n$. From here it can be easily seen that
\beq
\label{trPk}
\mathrm{trace}(P_k)=(n-k)s_k=c_kH_k,
\eeq
\beq
\label{trSPk}
\mathrm{trace}(S\circ P_k)=(k+1)s_{k+1=}c_kH_{k+1},
\eeq
and
\beq
\label{trS2Pk}
\mathrm{trace}(S^2\circ P_k)=(s_1s_{k+1}-(k+2)s_{k+2})=\bin{k+1}(nH_1H_{k+1}-(n-k-1)H_{k+2}),
\eeq
where
\[
c_k=(n-k)\bin{k}=(k+1)\bin{k+1}.
\]

These properties are all algebraic, and they can be found, for instance, in
\cite{Re}. There is still another non-algebraic property of $P_k$ that we need, which
can be found, for instance, in \cite[Lemma A]{Re1} and \cite[Equation (4.4)]{Ro} (see also \cite[page 118]{AG}).
The property we need is the following equation,
\beq
\label{trnabla}
\mathrm{tr}(P_k\circ \nabla_{X}S)=\g{\nabla s_{k+1}}{X}=\bin{k+1}\g{\nabla
H_{k+1}}{X}, \quad \mathrm{for} \quad  X\in\xm,
\eeq
where $\nabla S$ denotes the covariant differential of $S$,
\[
\nabla S(Y,X)=(\nabla_XS)Y=\nabla_X(SY)-S(\nabla_XY), \quad X,Y\in\xm.
\]

Associated to each Newton transformation $P_k$, we consider the second order linear
differential operator $L_k:\mathcal{C}^\infty(\m)\fle\mathcal{C}^\infty(\m)$ given by
\[
L_k(f)=\mathrm{trace}(P_k\circ\nabla^2f).
\]
Here $\nabla^2f:\mathcal{X}(\m)\fle\mathcal{X}(\m)$ denotes the
self-adjoint linear operator metrically equivalent to the hessian of $f$ and given by
\[
\g{\nabla^2f(X)}{Y}=\g{\nabla_X(\nabla f)}{Y}, \quad X,Y\in\xm.
\]
Consider $\{ E_1, \ldots, E_n \}$ a local orthonormal frame on \m\ and observe that
\begin{eqnarray*}
\mathrm{div}(P_k(\nabla f)) & = & \sum_{i=1}^n\g{(\nabla_{E_i}P_k)(\nabla f}{E_i}+
\sum_{i=1}^n\g{P_k(\nabla_{E_i}\nabla f)}{E_i}\\
{} & = & \g{\mathrm{div} P_k}{\nabla f}+L_k(f),
\end{eqnarray*}
where div denotes here the divergence on \m\ and
\[
\mathrm{div} P_k:=\mathrm{trace}(\nabla P_k)=\sum_{i=1}^n(\nabla_{E_i}P_k)(E_i).
\]
Obviously, $\mathrm{div}\ P_0=\mathrm{div}\ I=0$. Now Codazzi equation jointly with \rf{trnabla} imply that
$\mathrm{div} P_k=0$ also for every $k\geq 1$ \cite[Lemma B]{Re1}. To see it observe that, from the inductive
definition of $P_k$, we have
\[
(\nabla_{E_i}P_k)(E_i)=\bin{k}\g{\nabla H_k}{E_i}E_i-(\nabla_{E_i}S\circ P_{k-1})E_i-
(S\circ \nabla_{E_i}P_{k-1})E_i,
\]
so that
\[
\mathrm{div} P_k=\bin{k}\nabla H_k-\sum_{i=1}^n(\nabla_{E_i}S)(P_{k-1}E_i)-S(\mathrm{div} P_{k-1}).
\]
By Codazzi equation we know that $\nabla S$ is symmetric, and then for every $X\in\xm$
\begin{eqnarray*}
\sum_{i=1}^n\g{(\nabla_{E_i}S)(P_{k-1}E_i)}{X} & = & \sum_{i=1}^n\g{P_{k-1}E_i}{(\nabla_{E_i}S)X}
\sum_{i=1}^n\g{P_{k-1}E_i}{(\nabla_{X}S)E_i} \\
{} & = & \mathrm{tr}(P_{k-1}\circ \nabla_{X}S)=\bin{k}\g{\nabla H_k}{X}.
\end{eqnarray*}
In other words,
\[
\sum_{i=1}^n(\nabla_{E_i}S)(P_{k-1}E_i)=\bin{k}\nabla H_k,
\]
and then
\[
\mathrm{div} P_k=-S(\mathrm{div} P_{k-1}).
\]
Since $\mathrm{div} P_0=0$, this yields $\mathrm{div} P_k=0$ for every $k$. As a consequence,
$L_k(f)=\mathrm{div}(P_k(\nabla f))$ is a divergence form differential operator on \m.

\section{Examples}
\label{examples}
Let \x\ be an orientable hypersurface immersed into \mbar, with Gauss map $N$. For a fixed arbitrary vector $a\in\R{n+2}$,
let us consider the coordinate function $\g{a}{x}$ on \m. From $\nablao a=0$ we see that
\[
X(\g{a}{x})=\g{X}{a}=\g{X}{a^\top},
\]
for every vector field $X\in\xm$, where $a^\top\in\xm$ denotes the tangential
component of $a$,
\beq
\label{eq1}
a=a^\top+\g{a}{N}N+c\g{a}{x}x.
\eeq
Then the gradient of $\g{a}{x}$ on \m\ is given by $\nabla\g{a}{x}=a^\top$. By taking
covariant derivative in \rf{eq1} and using the Gauss and Weingarten formulae, we
also have from $\nablao a=0$ that
\beq
\label{eq2}
\nabla_X\nabla\g{a}{x}=\nabla_Xa^\top=\g{a}{N}SX-c\g{a}{x}X,
\eeq
for every tangent vector field $X\in\xm$. Therefore, by \rf{trSPk} we find that
\beq
\label{Lxa}
L_k\g{a}{x}=\g{a}{N}\mathrm{tr}(S\circ P_k)-c\g{a}{x}\mathrm{tr}(P_k)=
c_kH_{k+1}\g{a}{N}-cc_kH_k\g{a}{x}.
\eeq
That is
\beq
\label{Lx}
L_kx=c_kH_{k+1}N-cc_kH_kx.
\eeq

\begin{example}
\label{ex1}
It follows from \rf{Lx} that every hypersurface with vanishing $(k+1)$-th mean curvature and
having constant $k$-th mean curvature $H_k$ trivially satisfies \Lx\ with $A=-cc_k
H_k I_{n+2}\in\mathbb{R}^{(n+2)\times (n+2)}$ and $b=0$.
\end{example}

\begin{example}(Totally umbilical hypersurfaces in $\s{n+1}$).
\label{ex2}
As is well-known, the totally umbilical hypersurfaces of $\s{n+1}$ are the $n$-dimensional round spheres of radius
$0<r\leq 1$ which are obtained by intersecting $\s{n+1}$ with affine hyperplanes. Specifically, take $a\in\R{n+2}$ a unit
constant vector and, for a given $\tau\in(-1,1)$, let
\[
M_\tau=\{ x\in\s{n+1} : \g{a}{x}=\tau \}=\s{n}(\sqrt{1-\tau^2}).
\]
Then $M_\tau$ is a totally umbilical hypersurface in $\s{n+1}$ with Gauss map
$N(x)=(1/\sqrt{1-\tau^2})(a-\tau x)$ and shape operator $S=\tau/\sqrt{1-\tau^2}I$.
In particular, its higher order mean curvatures are given by
\[
H_{k}=\frac{\tau^k}{(1-\tau^2)^{k/2}}, \quad k=0,\ldots, n.
\]
Therefore, by equation \rf{Lx} we see that $M_\tau$ satisfies the condition \Lx\ for every $k=0,\ldots,n-1$, with
\[
A=\frac{-c_k\tau^k}{(1-\tau^2)^{(k+2)/2}}I_{n+2} \quad \mathrm{and} \quad
b=\frac{c_k\tau^{k+1}}{(1-\tau^2)^{(k+2)/2}}a.
\]
In particular, $b=0$ only when $\tau=0$, and then $M_0=\s{n}$ is a totally geodesic round sphere.
\end{example}

\begin{example}(Totally umbilical hypersurfaces in $\h{n+1}$).
\label{ex3}
Similarly to the case of the sphere, the totally umbilical hypersurfaces of $\h{n+1}$ are also obtained by intersecting
$\h{n+1}$ with affine hyperplanes of $\mathbb{R}^{n+2}_1$, but in this case there are three different types of hypersurfaces,
depending on the causal character of the hyperplane. To be more precise, take $a\in\mathbb{R}^{n+2}_1$ a non-zero constant
vector such that $\g{a}{a}\in\{1,0,-1\}$, and, for a given $\tau\in\mathbb{R}$, let
\[
M_\tau=\{ x\in\h{n+1} : \g{a}{x}=\tau \}.
\]
Then, when $\g{a}{a}+\tau^2>0$, $M_\tau$ is a totally umbilical hypersurface in $\h{n+1}$.
Observe that when $\g{a}{a}=1$ there is no restriction on the value of $\tau$ and $M_\tau=\h{n}(-\sqrt{1+\tau^2})$ is
a hyperbolic $n$-space of radius $-\sqrt{1+\tau^2}$. On the other hand, if $\g{a}{a}=-1$ then $|\tau|>1$ and
$M_\tau=\s{n}(\sqrt{\tau^2-1})$ is a round $n$-sphere of radius $\sqrt{\tau^2-1}$. Finally, when $\g{a}{a}=0$
then $\tau\neq 0$ and $M_\tau=\R{n}$ is a Euclidean space.

The Gauss map of $M_\tau$ is given by $N(x)=(1/\sqrt{\g{a}{a}+\tau^2})(a+\tau x)$, its shape operator is
$S=-\tau/\sqrt{\g{a}{a}+\tau^2}I$, and its higher order mean curvatures are given by
\[
H_{k}=\frac{(-1)^k\tau^k}{(\g{a}{a}+\tau^2)^{k/2}}, \quad k=0,\ldots, n.
\]
Therefore, by equation \rf{Lx} we see that $M_\tau$ satisfies the condition \Lx\ for every $k=0,\ldots,n-1$, with
\[
A=\frac{(-1)^{k}c_k\g{a}{a}\tau^k}{(\g{a}{a}+\tau^2)^{(k+2)/2}}I_{n+2} \quad \mathrm{and} \quad
b=\frac{(-1)^{k+1}c_k\tau^{k+1}}{(\g{a}{a}+\tau^2)^{(k+2)/2}}a.
\]
In particular, $b=0$ only when $\tau=0$, and then $M_0=\h{n}$ is a totally geodesic hyperbolic space. On the other
hand, the totally umbilical Euclidean spaces in \h{n+1} (corresponding to the case $\g{a}{a}=0$) satisfy
the condition \Lx\ with $A=0$.
\end{example}

\begin{example}(Standard Riemannian products in $\s{n+1}$ and $\h{n+1}$)
\label{ex4}
Here we will consider the case where \m\ is a standard Riemannian product; that is, \m\ is either the Riemannian
product $\s{m}(\sqrt{1-r^2})\times\s{n-m}(r)\subset\s{n+1}$ with $0<r<1$, or the Riemannian product
$\h{m}(-\sqrt{1+r^2})\times\s{n-m}(r)\subset\h{n+1}$ with $r>0$, for a certain $m=1,\ldots,n-1$. After a rigid motion of
the ambient space, we may consider that \m\ is defined by the equation
\[
M=\{ x\in\mbar\ : x_{m+1}^2+\cdots x_{n+1}^2=r^2 \}.
\]
In that case, the Gauss map on \m\ is
\[
N(x)=\left(\frac{-cr}{\sqrt{1-cr^2}}x_0,\ldots,\frac{-cr}{\sqrt{1-cr^2}}x_m,\frac{\sqrt{1-cr^2}}{r}x_{m+1},\ldots,\frac{\sqrt{1-cr^2}}{r}x_{n+1}\right).
\]
and its the principal curvatures are
\[
\kappa_1=\cdots=\kappa_m=\frac{cr}{\sqrt{1-cr^2}}, \quad \kappa_{m+1}=\cdots=\kappa_n=\frac{-\sqrt{1-cr^2}}{r}.
\]
In particular, the higher order mean curvatures are all constant.
Therefore, using \rf{Lx} we get that
\[
L_kx=
(\lambda x_0,\ldots,\lambda x_m,\mu x_{m+1},\ldots,\mu x_{n+1})
\]
where $\lambda$ and $\mu$ are both constants,
\[
\lambda=\frac{-cc_kH_{k+1}r}{\sqrt{1-cr^2}}-cc_kH_k, \quad
\mu=\frac{c_kH_{k+1}\sqrt{1-cr^2}}{r}-cc_kH_k.
\]
That is, \m\ satisfies the condition \Lx\ with $b=0$ and
\[
A=\mathrm{diag}[\lambda,\ldots,\lambda,\mu,\ldots,\mu].
\]
\end{example}

\section{Some computations and first auxiliary results}
In Section\rl{examples} we have computed the operator $L_k$ acting on the coordinate functions of a hypersurface.
On the other hand, consider now the coordinate functions of its Gauss map $N$, that is, the function $\g{a}{N}$ on \m,
where $a\in\R{n+2}$ is a fixed arbitrary vector. From $\nablao a=0$ we also see that
\[
X(\g{a}{N})=-\g{SX}{a}=-\g{X}{S(a^\top)}
\]
for every vector field $X\in\xm$, so that
\[
\nabla\g{a}{N}=-S(a^\top).
\]
Therefore, from \rf{eq2} we get
\begin{eqnarray}
\label{eq3}
\nonumber \nabla_X(\nabla\g{a}{N}) & = & -\nabla_X(Sa^\top)=-\nabla S(a^\top,X)-S(\nabla_Xa^\top)\\
{} & = & -(\nabla_XS)a^\top-\g{a}{N}S^2X+c\g{a}{x}SX.
\end{eqnarray}
By Codazzi equation we know that $\nabla S$ is symmetric and then
\[
\nabla S(a^\top,X)=\nabla S(X,a^\top)=(\nabla_{a^\top}S)X.
\]
Therefore using this in \rf{eq3}, jointly with \rf{trS2Pk} and \rf{trnabla}, we get
\begin{eqnarray}
\label{LNa}
\nonumber L_k\g{a}{N} & = &
-\mathrm{tr}(P_k\circ \nabla_{a^\top}S)-\g{a}{N}\mathrm{tr}(S^2\circ P_k)+c\g{a}{x}\mathrm{tr}(S\circ P_k) \\
{} & = & -\bin{k+1}\g{\nabla H_{k+1}}{a}\\
\nonumber {} & {} & -\bin{k+1}(nH_1H_{k+1}-(n-k-1)H_{k+2})\g{a}{N}\\
\nonumber {} & {} & +cc_kH_{k+1}\g{a}{x}.
\end{eqnarray}
In other words,
\begin{eqnarray}
\label{LN}
\nonumber L_kN & = & -\bin{k+1}\nabla H_{k+1} \\
{} & - & \bin{k+1}(nH_1H_{k+1}-(n-k-1)H_{k+2})N \\
\nonumber {} & + & \bin{k+1}c(k+1)H_{k+1}x.
\end{eqnarray}

Let us assume that, for a fixed $k=0,\ldots,n-1$, the immersion
$$
\x\
$$
satisfies the condition
\beq
\label{eq4}
\Lx,
\eeq
for a constant matrix $A\in\R{(n+2)\times(n+2)}$ and a constant vector $b\in\R{n+2}$. From \rf{Lx} we get that
\beq
\label{Ax}
Ax=-b+c_kH_{k+1}N-cc_kH_{k}x=-b^\top+(c_kH_{k+1}-\g{b}{N})N-c(c_kH_k+\g{b}{x})x,
\eeq
where $b^\top\in\xm$ denotes the tangential component of $b$. Now, if we take covariant derivative
in \rf{eq4} and use the equation \rf{Lx} as well as Weingarten formula, we obtain
\beq
\label{AX}
AX=-c_kH_{k+1}SX-cc_kH_{k}X+c_k\g{\nabla H_{k+1}}{X}N-cc_k\g{\nabla H_{k}}{X}x
\eeq
for every tangent vector field $X\in\xm$. On the other hand, taking into account that
\[
L_k(fg)=(L_kf)g+f(L_kg)+2\g{P_k(\nabla f)}{\nabla g}, \quad f,g\in\mathcal{C}^\infty(\m),
\]
we also get from \rf{Lxa} and \rf{LNa} that
\begin{eqnarray*}
L_k(L_k\g{a}{x})=
-c_k\bin{k+1}H_{k+1}\g{\nabla H_{k+1}}{a}-2c_k\g{(S\circ P_k)(\nabla H_{k+1})}{a}-2cc_k\g{P_k(\nabla H_{k})}{a}\\
-c_k\left(\bin{k+1}H_{k+1}(nH_1H_{k+1}-(n-k-1)H_{k+2})+cc_kH_kH_{k+1}-L_kH_{k+1}\right)\g{a}{N}\\
+c_k\left(cc_kH^2_{k+1}+c_kH^2_k-cL_kH_{k}\right)\g{a}{x}.
\end{eqnarray*}
Equivalently,
\begin{eqnarray*}
L_k(L_kx) & = &
-c_k\bin{k+1}H_{k+1}\nabla H_{k+1}-2c_k(S\circ P_k)(\nabla H_{k+1})-2cc_kP_k(\nabla H_{k})\\
{} & {} & -c_k\left(\bin{k+1}H_{k+1}(nH_1H_{k+1}-(n-k-1)H_{k+2})+cc_kH_kH_{k+1}-L_kH_{k+1}\right)N\\
{} & {} & +c_k\left(cc_kH^2_{k+1}+c_kH^2_k-cL_kH_{k}\right)x.
\end{eqnarray*}
From here, by applying the operator $L_k$ on both sides of \rf{eq4} and using again \rf{Lx}, we
have
\begin{eqnarray}
\label{AN}
\nonumber H_{k+1}AN & = &
-\bin{k+1}H_{k+1}\nabla H_{k+1}-2(S\circ P_k)(\nabla H_{k+1})-2cP_k(\nabla H_{k}) \\
{} & {} & -\left(\bin{k+1}H_{k+1}(nH_1H_{k+1}-(n-k-1)H_{k+2})+cc_kH_kH_{k+1}-L_kH_{k+1}\right)N\\
\nonumber {} & {} & +\left(cc_kH^2_{k+1}+c_kH^2_k-cL_kH_{k}\right)x+cH_kAx.
\end{eqnarray}
Using here \rf{Ax}, we get
\begin{eqnarray}
\label{ANbis}
\nonumber H_{k+1}AN & = &
-\bin{k+1}H_{k+1}\nabla H_{k+1}-2(S\circ P_k)(\nabla H_{k+1})-2cP_k(\nabla H_{k})-cH_kb^\top \\
{} & {} & -\left(\bin{k+1}H_{k+1}(nH_1H_{k+1}-(n-k-1)H_{k+2})+cH_k\g{b}{N}-L_kH_{k+1}\right)N\\
\nonumber {} & {} & +\left(cc_kH^2_{k+1}-cH_k\g{b}{x}-cL_kH_{k}\right)x.
\end{eqnarray}

\subsection{The case where $A$ is self-adjoint}
From \rf{AX} we have
\beq
\label{eq5bis}
\g{AX}{Y}=\g{X}{AY}
\eeq
for every tangent vector fields $X,Y\in\xm$. In other words, the endomorphism determined by $A$ is always self-adjoint
when restricted to the tangent hyperplanes of the hypersurface. Therefore, $A$ is self-adjoint if and only if the three
following equalities hold
\beq
\label{EQ1}
\g{AX}{x}=\g{x}{AX} \quad \mbox{for every $X\in\xm$},
\eeq
\beq
\label{EQ2}
\g{AX}{N}=\g{X}{AN} \quad \mbox{for every $X\in\xm$},
\eeq
and
\beq
\label{EQ3}
\g{AN}{x}=\g{N}{Ax}.
\eeq

From \rf{Ax} and \rf{AX} it easily follows that \rf{EQ1} is equivalent to
\beq
\label{EQ1bis}
\nabla\g{b}{x}=b^\top=c_k\nabla H_k,
\eeq
that is, $\g{b}{x}-c_kH_k$ is constant on \m. On the other hand, from \rf{AX} and \rf{AN}, and using also \rf{EQ1bis}, it
follows that, at points where $H_{k+1}\neq 0$, \rf{EQ2} is equivalent to
\begin{eqnarray}
\label{EQ2bis}
\frac{2}{H_{k+1}}(S\circ P_k)(\nabla H_{k+1})+(k+2)\bin{k+1}\nabla H_{k+1}=\\
\nonumber -\frac{c}{H_{k+1}}\left(2P_k(\nabla H_{k})+c_kH_k\nabla H_{k}\right).
\end{eqnarray}
Finally, using again \rf{EQ1bis} we have by \rf{Lxa} that
\beq
\label{LkHk}
L_kH_k=\frac{1}{c_k}L_k\g{b}{x}=H_{k+1}\g{b}{N}-cH_k\g{b}{x}.
\eeq
Observe also that
\[
\g{Ax}{x}=-\g{b}{x}-cc_kH_k.
\]
Therefore, from \rf{AN} we get that
\begin{eqnarray*}
H_{k+1}\g{AN}{x} & = & c_kH_{k+1}^2+cc_kH_k^2-L_kH_k+cH_k\g{Ax}{x}\\
{} & = & c_kH_{k+1}^2-H_{k+1}\g{b}{N}\\
{} & = & H_{k+1}\g{N}{Ax}.
\end{eqnarray*}
Thus we have that, at points where $H_{k+1}\neq 0$, the first two equalities \rf{EQ1} and \rf{EQ2} imply the third one
\rf{EQ3}.

Now we are ready to prove the following auxiliary result.
\begin{lemma}
\label{CMCbiss}
Let \x\ be an orientable hypersurface satisfying the condition \Lx,
for some self-adjoint constant matrix $A\in\R{(n+2)\times(n+2)}$ and some constant vector $b\in\R{n+2}$. Then $H_{k}$
is constant if and only if $H_{k+1}$ is constant.
\end{lemma}
\begin{proof}
Assume that $H_k$ is constant and let us consider the open set
\[
\mathcal{U}=\{ p\in\m : \nabla H_{k+1}^2(p)\neq 0 \}.
\]
Our objective is to show that $\mathcal{U}$ is empty. Assume that $\mathcal{U}$ is non-empty. From \rf{EQ2bis} we have
that
\[
\frac{2}{H_{k+1}}(S\circ P_k)(\nabla H_{k+1})+(k+2)\bin{k+1}\nabla H_{k+1}=0  \quad \mathrm{on} \quad \mathcal{U}.
\]
Equivalently,
\[
(S\circ P_k)(\nabla H_{k+1})=-\frac{k+2}{2}\bin{k+1}H_{k+1}\nabla H_{k+1} \quad \mathrm{on} \quad \mathcal{U}.
\]
Then, reasoning exactly as Al\'\i as and G\"urb\"uz in \cite[Lemma 5]{AG} (starting from equation (23) in \cite{AG}) we conclude
that $H_{k+1}$ is locally constant on $\mathcal{U}$, which is a contradiction. Actually, the proof in \cite{AG} works also
here word by word, with the only observation that, since we are assuming that $H_k$ is constant, \rf{AX} reduces now to
\[
AX=-c_kH_{k+1}SX-cc_kH_{k}X+c_k\g{\nabla H_{k+1}}{X}N.
\]
Therefore, instead of having $AE_i=-c_kH_{k+1}\kappa_iE_i$, now we have
$AE_i=-c_k(H_{k+1}\kappa_i+cH_k)E_i$ for every $m+1\leq i\leq n$ (see the last paragraph of the proof of
\cite[Lemma 5]{AG}). But $H_k$ being constant, that makes no difference to the reasoning.

Conversely, assume now that $H_{k+1}$ is constant and let us consider the open set
\[
\mathcal{V}=\{ p\in\m : \nabla H_{k}^2(p)\neq 0 \}.
\]
Our objective now is to show that $\mathcal{V}$ is empty. Let us consider first the case where $H_{k+1}=0$ and
assume that $\mathcal{V}$ is non-empty. In this case, by \rf{EQ1bis} and \rf{LkHk}, \rf{ANbis} reduces to
\[
-2cP_k(\nabla H_k)-cc_kH_k\nabla H_k-cH_k\g{b}{N}N=0.
\]
Thus, $\g{b}{N}=0$ on $\mathcal{V}$. By \rf{Ax} this gives $\g{AN}{x}=\g{N}{Ax}=0$ and, since
$\g{AN}{X}=\g{N}{AX}=0$ for every $X\in\xm$, we obtain that $AN=\g{AN}{N}N$; that is, $N$ is an eigenvector of $A$
with corresponding eigenvalue $\lambda=\g{AN}{N}$. In particular, $\lambda$ is locally constant on $\mathcal{V}$.
Therefore,
\begin{eqnarray*}
AX & = & -cc_kH_{k}X-cc_k\g{\nabla H_{k}}{X}x \\
AN & = & \lambda N \\
Ax & = & -c_k\nabla H_k -c(2c_kH_k+\alpha)x,
\end{eqnarray*}
where $\alpha=\g{b}{x}-c_kH_k$ and $\lambda$ are both locally constant on $\mathcal{V}$. Then,
\[
\mathrm{tr}(A)=-ncc_kH_k+\lambda-c(2c_kH_k+\alpha)=\mathrm{constant},
\]
which implies that $H_k$ is locally constant on $\mathcal{V}$, which is a contradiction.

On the other hand, if $H_{k+1}\neq 0$ is constant and we assume that $\mathcal{V}$ is non-empty, then from \rf{EQ2bis} we
have that
\[
2P_k(\nabla H_{k})+c_kH_k\nabla H_{k}=0 \quad \mathrm{on} \quad \mathcal{V}.
\]
Equivalently,
\beq
\label{eq10}
P_k(\nabla H_{k})=-\frac{c_k}{2}H_{k}\nabla H_{k} \quad \mathrm{on} \quad \mathcal{V}.
\eeq
Here, we will follow a similar reasoning to that in \cite[Lemma 5]{AG}. Consider
$\{ E_1,\ldots, E_n\}$ a local orthonormal frame of principal directions of $S$ such that $SE_i=\kappa_iE_i$ for
every $i=1,\ldots,n$, and then
\[
P_{k}E_i=\mu_{i,k}E_i,
\]
with
\beq
\label{eq9.5}
\mu_{i,k}=\sum_{j=0}^{k}(-1)^{j}\bin{k-j}H_{k-j}\kappa_i^{j}=
\sum_{i_1<\cdots<i_{k},i_j\neq i}\kappa_{i_1}\cdots\kappa_{i_{k}}.
\eeq
Therefore, writing
\[
\nabla H_{k}=\sum_{i=1}^{n}\g{\nabla H_{k}}{E_i}E_i
\]
we see that \rf{eq10} is equivalent to
\[
\g{\nabla H_{k}}{E_i}\left(\mu_{i,k}+\frac{c_k}{2}H_{k}\right)=0 \quad \mathrm{on} \quad \mathcal{V}
\]
for every $i=1,\ldots,n$. Thus, for every $i$ such that $\g{\nabla H_{k}}{E_i}\neq 0$ on $\mathcal{V}$ we get
\beq
\label{eq11}
\mu_{i,k}=-\frac{c_k}{2}H_{k}.
\eeq
This implies that $\g{\nabla H_{k}}{E_i}=0$ necessarily for some $i$. Otherwise, we would have \rf{eq11} for every
$i=1,\ldots,n$, which would imply
\[
c_kH_{k}=\mathrm{tr}(P_{k})=\sum_{i=1}^{n}\mu_{i,k}=-\frac{nc_k}{2}H_{k},
\]
and thus $H_{k}=0$ on $\mathcal{V}$, which is a contradiction.

Therefore, re-arranging the local orthonormal frame if necessary, we may assume that for some $1\leq m<n$ we have
$\g{\nabla H_{k}}{E_i}\neq 0$ for $i=1,\ldots,m$, $\g{\nabla H_{k}}{E_i}=0$ for $i=m+1,\ldots,n$, and
$\kappa_1<\kappa_2<\cdots<\kappa_m$. The integer $m$ measures the number of linearly independent principal
directions of $\nabla H_{k}$, and $\nabla H_{k}$ is a principal direction of $S$ if and only if $m=1$.
From \rf{eq11} we know that
\beq
\label{eq12}
\mu_{1,k}=\cdots=\mu_{m,k}=-\frac{c_k}{2}H_{k}\neq 0 \quad \mathrm{on} \quad \mathcal{V}.
\eeq
Thus, by \rf{eq9.5} it follows that $\kappa_1<\kappa_2<\cdots<\kappa_m$ are $m$ distinct real roots of the
following polynomial equation of degree $k$,
\[
Q(t)=\sum_{j=0}^{k}(-1)^{j}\bin{k-j}H_{k-j}t^{j}=-\frac{c_k}{2}H_{k}.
\]
In particular $m\leq k$. On the other hand, each $\kappa_i$ is also a root of the characteristic polynomial of
$S$, which can be written as
\[
Q_S(t)=(-1)^{k}t^{n-k}Q(t)+\sum_{j=k+1}^{n}(-1)^{j}\bin{j}H_{j}t^{n-j}.
\]
Then, $\kappa_1<\kappa_2<\cdots<\kappa_m$ are also $m$ distinct real roots of the
following polynomial equation of degree $n-k$,
\[
(-1)^{k+1}\frac{c_k}{2}H_{k}t^{n-k}+\sum_{j=k+1}^{n}(-1)^{j}\bin{j}H_{j}t^{n-j}=0.
\]
In particular, $m\leq n-k$, that is, $n-m\geq k$. Now we claim that
\beq
\label{induccion}
\mu_{1,k}=\cdots=\mu_{m,k}=\sum_{m<i_1<\cdots<i_{k}}\kappa_{i_1}\cdots\kappa_{i_{k}}.
\eeq
The proof of \rf{induccion} follows exactly as the proof of equation (29) in \cite{AG} and we omit it here.

Finally, from equation \rf{AX} we have
\[
AE_i=-c_k(H_{k+1}\kappa_i+cH_k)E_i
\]
for every $m+1\leq i\leq n$. Therefore, every $-c_k(H_{k+1}\kappa_i+cH_k)$ with $i=m+1,\ldots n$ is a constant eigenvalue
$\alpha_i$ of the constant matrix $A$. Then,
\[
\kappa_i=-\frac{\alpha_i+cc_kH_k}{c_kH_{k+1}} \quad \mbox{for every $i=m+1,\ldots n$}
\]
and from \rf{induccion} and \rf{eq12} we get that
\[
-\frac{c_k}{2}H_{k}=\sum_{m<i_1<\cdots<i_{k}}\kappa_{i_1}\cdots\kappa_{i_{k}}=
\frac{(-1)^{k}}{c_k^{k}H_{k+1}^{k}}\sum_{m<i_1<\cdots<i_{k}}(\alpha_{i_1}+cc_kH_k)\cdots(\alpha_{i_{k}}+cc_kH_k)
\]
on $\mathcal{V}$. But this means that $H_{k}$ is locally constant on $\mathcal{V}$, which is a contradiction with the
definition of $\mathcal{V}$. This finishes the proof of Lemma\rl{CMCbiss}.
\end{proof}

\section{Proof of Theorem\rl{th1}}
We have already checked in Section\rl{examples} that each one of the
hypersurfaces mentioned in Theorem\rl{th1} does satisfy the condition \Lxb\ for a self-adjoint
constant matrix $A$. Conversely, let us
assume that \x\ satisfies the condition \Lxb\ for some self-adjoint constant matrix $A\in\R{(n+2)\times(n+2)}$.
Since $b=0$, from \rf{EQ1bis} we get that $H_k$ is constant on \m. Thus, by Lemma\rl{CMCbiss} we know that $H_{k+1}$ is
also constant on \m. If $H_{k+1}=0$ there is nothing to prove. Then, we may assume that $H_{k+1}$ is a non-zero constant
and $H_k$ is also constant. Then from \rf{AX} and \rf{AN} we obtain
\beq
\label{AX1}
AX=-c_kH_{k+1}SX-cc_kH_kX
\eeq
for every tangent vector field $X\in\xm$, and
\beq
\label{AN1}
AN=\alpha N+c_k\left(cH_{k+1}+\frac{H^2_k}{H_{k+1}}\right)x+c\frac{H_k}{H_{k+1}}Ax,
\eeq
with
\[
\alpha=-\bin{k+1}(nH_1H_{k+1}-(n-k-1)H_{k+2})-cc_kH_k.
\]
Taking covariant derivative in \rf{AN1} and using \rf{AX1} we have for every $X\in\xm$
\begin{eqnarray*}
\nablao_X(AN) & = & \g{\nabla\alpha}{X}N-\alpha SX+c_k\left(cH_{k+1}+\frac{H^2_k}{H_{k+1}}\right)X+
c\frac{H_k}{H_{k+1}}AX\\
{} & = & \g{\nabla\alpha}{X}N+\bin{k+1}(nH_1H_{k+1}-(n-k-1)H_{k+2} SX+cc_kH_{k+1}X.
\end{eqnarray*}
On the other hand, from \rf{AX1} we also find that
\[
\nablao_X(AN)=A(\nablao_XN)=-A(SX)=c_kH_{k+1}S^2X+cc_kH_kSX
\]
It follows from here that $\g{\nabla\alpha}{X}=0$ for every $X\in\xm$, that is, $\alpha$ is constant on \m,
and also that the shape operator $S$ satisfies the following quadratic equation
\[
S^2+\lambda S-cI=0,
\]
where
\[
\lambda=\frac{\alpha}{c_kH_{k+1}}+2c\frac{H_k}{H_{k+1}}=\mathrm{constant}.
\]
As a consequence, either \m\ is totally umbilical in \mbar\ (but not totally geodesic, because of $H_{k+1}\neq 0$) or \m\
is an isoparametric hypersurface of
\mbar\ with two constant principal curvatures. The former cannot occur, because the only totally umbilical hypersurfaces in
\mbar\ which satisfy \Lxb\ with $b=0$ are the totally geodesic ones (see Examples\rl{ex2} and\rl{ex3}). In the latter, from
well-known results by Lawson
\cite[Lemma 2]{La} and Ryan \cite[Theorem 2.5]{Ry} we conclude that \m\ is an open piece of a standard Riemannian product.

\section{Proof of Theorem\rl{th2}}
We have already checked in Section\rl{examples} that each one of the
hypersurfaces mentioned in Theorem\rl{th2} does satisfy the condition \Lx\ for a self-adjoint
constant matrix $A$. Conversely, let us assume
that \x\ satisfies the condition \Lx\ for some self-adjoint constant matrix $A\in\R{(n+2)\times(n+2)}$ and some
non-zero constant vector $b\in\R{n+2}$. Since $H_k$ is assumed to be constant, by Lemma\rl{CMCbiss} we know that $H_k$ and
$H_{k+1}$ are both constant on \m. The case $H_{k+1}=0$ cannot occur, because in that case we have $b=0$
(Example\rl{ex1}). Therefore, we have that $H_{k+1}$ is a non-zero constant and $H_k$ is also constant.
Then from \rf{AX} and \rf{AN} we obtain
\beq
\label{AX1b}
AX=-c_kH_{k+1}SX-cc_kH_kX
\eeq
for every tangent vector field $X\in\xm$, and
\beq
\label{AN1b}
AN=\alpha N+c_k\left(cH_{k+1}+\frac{H^2_k}{H_{k+1}}\right)x+c\frac{H_k}{H_{k+1}}Ax,
\eeq
with
\[
\alpha=-\bin{k+1}(nH_1H_{k+1}-(n-k-1)H_{k+2})-cc_kH_k.
\]
Taking covariant derivative in \rf{AN1b} and using \rf{AX1b} we have for every $X\in\xm$
\begin{eqnarray*}
\nablao_X(AN) & = & \g{\nabla\alpha}{X}N-\alpha SX+c_k\left(cH_{k+1}+\frac{H^2_k}{H_{k+1}}\right)X+
c\frac{H_k}{H_{k+1}}AX\\
{} & = & \g{\nabla\alpha}{X}N+\bin{k+1}(nH_1H_{k+1}-(n-k-1)H_{k+2} SX+cc_kH_{k+1}X.
\end{eqnarray*}
On the other hand, from \rf{AX1b} we also find that
\[
\nablao_X(AN)=A(\nablao_XN)=-A(SX)=c_kH_{k+1}S^2X+cc_kH_kSX
\]
It follows from here that $\g{\nabla\alpha}{X}=0$ for every $X\in\xm$, that is, $\alpha$ is constant on \m,
and also that the shape operator $S$ satisfies the following quadratic equation
\[
S^2+\lambda S-cI=0,
\]
where
\[
\lambda=\frac{\alpha}{c_kH_{k+1}}+2c\frac{H_k}{H_{k+1}}=\mathrm{constant}.
\]
As a consequence, either \m\ is totally umbilical in \mbar\ or \m\ is an isoparametric hypersurface of
\mbar\ with two constant principal curvatures. In the latter, from well-known results by Lawson \cite[Lemma 2]{La} and
Ryan \cite[Theorem 2.5]{Ry} we would get that \m\ is an open piece of a standard Riemannian product, but this case
cannot occur because they satisfy the condition \Lx\ with $b=0$ (Example\rl{ex4}).

\section*{Acknowledgements}
This work was done while the second author was spending his sabbatical leave at the Warwick Mathematics Institute (WMI).
He wants to thank Tarbiat Modarres University for its financial support and the WMI for its hospitality. He also would
like to thank the first author for his visit to WMI for invaluable discussion about this work. The
authors thank to the referee for valuable suggestions which improved the paper.

\bibliographystyle{amsplain}

\end{document}